\documentclass[12pt]{amsart}
\usepackage{amsmath}
\usepackage{footnote}
\usepackage{amssymb}

\newcommand{\tr}{\text{Tr}}

\newcommand{\Z}{{\mathbb Z}}
\newcommand{\Q}{{\mathbb Q}}

\newcommand{\F}{{\mathbb F}}

\newcommand{\Tr}{{\mathrm{Tr}}}
\newcommand{\Norm}{{\mathrm{Norm}}}
\newcommand{\Cay}{{\mathrm{Cay}}}

\newcommand{\PG}{\mathsf{PG}}
\newcommand{\Aut}{\mathsf{Aut}}
\newcommand{\Ha}{\mathsf{H}}
\newcommand{\Wa}{\mathsf{W}}
\renewcommand{\mod}{{\mathrm{mod\, \, }}}

\newcommand{\cQ}{\mathsf{Q}}

\theoremstyle{plain}
\newtheorem{theorem}{Theorem}[section]

\newtheorem{proposition}[theorem]{Proposition}
\newtheorem{example}[theorem]{Example}
\newtheorem{conj}[theorem]{Conjecture}

\numberwithin{equation}{section}

\theoremstyle{remark}
\newtheorem{remark}[theorem]{Remark}

\usepackage{url, hyperref}

\renewcommand\le{\leqslant}
\renewcommand\ge{\geqslant}

\def\qi#1 {\fbox {\footnote {\ }}\ \footnotetext { From Qi: {\color{blue}#1}}}

\begin{document}

\title[Cyclotomy, difference sets, sequences, and srgs]{Cyclotomy,  difference sets, sequences with low correlation, strongly regular graphs, and related geometric substructures}

\author{Koji Momihara}
\address{Division of Natural Science, Faculty of Advanced Science and Technology, Kumamoto University, 2-40-1 Kurokami, Kumamoto 860-8555, Japan}
\curraddr{}
\email{momihara@educ.kumamoto-u.ac.jp}
\thanks{}

\author{Qi Wang}
\address{Department of Computer Science and Engineering, 
Southern University of Science and Technology, 
Shenzhen, Guangdong 518055, China}
\curraddr{}
\email{wangqi@sustc.edu.cn}
\thanks{}

\author{Qing Xiang}
\address{Department of Mathematical Sciences, 
University of Delaware, 
Newark DE 19716, USA}
\curraddr{}
\email{qxiang@udel.edu}
\thanks{Koji Momihara was supported by 
JSPS under Grant-in-Aid for Young Scientists (B) 17K14236 and Scientific Research (B) 15H03636. Qi Wang was supported by the National Natural Science Foundation of China under Grant no. 11601220. Qing Xiang was supported by an NSF grant DMS-1600850, and a JSPS invitational fellowship for research in Japan S17114.}



\date{}

\dedicatory{}






\begin{abstract}
  In this paper, we survey constructions of and nonexistence results on combinatorial/geometric structures which arise from unions of cyclotomic classes of finite fields. In particular, we survey both classical and recent results on difference sets related to cyclotomy, and cyclotomic constructions of sequences with low correlation.  We also give an extensive survey of recent results on constructions of strongly regular Cayley graphs and related geometric substructures such as $m$-ovoids and $i$-tight sets in classical polar spaces.
\end{abstract}

\maketitle










\section{Introduction}\label{sec-intro}

Let $q = p^\ell$ be a prime power, and $\F_q$ be the finite field of order $q$. We use $\F_q^*$ to denote the set of nonzero elements of $\F_q$. It is well known that $\F_q^*$ is a cyclic group of order $q-1$. When $q$ is odd, let $C_0$ denote the unique subgroup of index 2 of $\F_q^*$; that is, $C_0$ is the subgroup of $\F_q^*$ consisting of the nonzero squares of $\F_q$. The set $C_0$ has played very important roles in the construction of various combinatorial structures such
as Hadamard matrices, difference sets, and strongly regular graphs. The earliest use of $C_0$ for constructing Hadamard matrices goes back to Paley~\cite{Paley33}. Subsequently, many researchers considered using subgroups of $\F_q^*$ of higher indices and their cosets for constructing difference sets, binary sequences with low correlation, and strongly regular Cayley graphs, etc. The additive properties of the subgroups of $\F_q^*$ form a large part of what we call the
theory of cyclotomy today. To a large extent, the theory of cyclotomy is a study of generalizations of Paley's work in~\cite{Paley33}.

We now give the definition of difference sets in a (not necessarily cyclic) group of order $v$. Let $G$ be a finite multiplicative group of order $v$. A $k$-element subset $D$ of $G$ is called a {\em $(v,k,\lambda)$ difference set\/} in $G$ if the list of ``differences'' $d_1d_2^{-1}$, $d_1,d_2\in D$, $d_1\neq d_2$, represents each nonidentity element in $G$ exactly $\lambda$ times. A moment's reflection shows that the translates of $D$ by all group elements form the blocks of a $(v,k,\lambda)$ symmetric design, and $G$ is a regular automorphism group of the design. For this reason difference sets play an important role in combinatorial design theory. 

Given a subset $D$ in the cyclic group $(\Z/v\Z, +)$, we define its {\em characteristic sequence} ${\bf s}=(s_i)_{0\leq i\leq v-1}$ with the {\em support} $D$ by setting $s_i=1$ if $i\in D$, and $s_i=-1$ otherwise. The {\em periodic autocorrelation} of a binary sequence ${\bf s}$ at the shift $\tau$, $0 \leq \tau < v$, is defined as ${\mathcal A}_{\bf s} (\tau) = \sum_{i=0}^{v-1} s_i s_{i+\tau}$, where $i + \tau$ is read modulo the period $v$. From the definition of difference set, we see that $D$ is a $(v,k,\lambda)$ difference set in $\Z/v\Z$ if and only if 
\begin{equation}\label{eqn-cs} 
  {\mathcal A}_{\bf s}(\tau)= \left\{
\begin{array}{ll}
  v, & \mbox{if} \; \tau\equiv  0 \pmod v, \\
v-4(k-\lambda), & \mbox{otherwise}. \end{array} \right. 
\end{equation}

This shows the equivalence of binary sequences with two-level autocorrelation and cyclic $(v,k,\lambda)$ difference sets. More generally, $(v,k,\lambda)$ abelian difference sets are equivalent to binary arrays with two-level autocorrelation. For background material on difference sets, we refer the reader to the books~\cite{Bau71, Lan83} and Chapter 6 of~\cite{BJL99}.

Let $q=p^{\ell}$ be a prime power, and let $\gamma$ be a fixed primitive element of $\F_q$. Let $N>1$ be a divisor of $q-1$. We define the $N^{\rm th}$ {\em cyclotomic classes} 
$C_i^{(N,q)}$ of $\F_q$ by
$$C_i^{(N,q)}=\{\gamma^{jN+i}\mid 0\leq j\leq \frac {q-1}{N}-1\},$$ where $0\leq i\leq N-1$. That is, $C_0^{(N,q)}$ is the subgroup of $\F_q^*$ consisting of all nonzero $N^{\rm th}$ powers in $\F_q$, and $C_i^{(N,q)}=\gamma^i C_0^{(N,q)}$, for $1\leq i\leq N-1$. The case where $N=2$ was first used by Paley~\cite{Paley33} to construct the Paley difference set when $q\equiv 3\pmod 4$, and the Paley graph when $q\equiv 1\pmod 4$. Even though the construction is
deterministic, the resulting combinatorial structures (i.e., the Paley difference sets/graphs) are pseudorandom or quasirandom. The $N^{\rm th}$ cyclotomic classes (with $N>2$) also exhibit pseudorandom behaviors. 

\begin{itemize}
  \item[(1)] Roughly speaking, a pseudorandom graph is a graph that behaves like a random graph of the same edge density. The notion of quasirandom (also called pseudorandom) graphs was made precise by Thomason \cite{Tho87} and Chung, Graham and Wilson \cite{CGW89}. The Paley graphs are now standard examples of explicitly constructed quasirandom graphs.
  \item[(2)] Elements of $C_0^{(N,q)}$ are distributed in $\F_q$ in a way that is random-like and also very regular at the same time. Here by random-like behavior, we mean that "being an $N^{\rm th}$ power" is like a random event of probability $\frac{1}{N}$. For the precise statement we refer the reader to Sziklai~\cite{Szi01} (see also~\cite{Win98}). The $N = 2$ case was treated by Sz\"onyi~\cite{Szo92} and Babai, Gal and Widgerson~\cite{BGW99}.  
  \item[(3)] The characteristic sequences of many difference sets from cyclotomic classes are pseudorandom with respect to certain randomness postulates, including balancedness, run property, low autocorrelation~\cite{Gol82}, pattern distribution~\cite{Ding98}, etc.   
\end{itemize}
In this survey paper, we will mainly focus on constructions of various combinatorial/geometric structures by using cyclotomic classes. The paper is organized as follows. In Section~\ref{cycdiff}, we survey both classical and recent results on difference sets related to cyclotomy. The highlights are some recent results of Xia \cite{Xia18} on the long-standing conjecture that if $C_0^{(N,q)}$ is a difference set in $(\F_q, +)$, then $N$ is a power of $2$; and the
constructions of skew Hadamard difference sets by Feng and the third author \cite{FX12} by using unions of cyclotomic classes. In Section 3, we give a brief survey of results on sequences with low correlation which are related to cyclotomy. Section 4 is devoted to strongly regular Cayley graphs arising from cyclotomy and related geometric substructures such as $m$-ovoids and $i$-tight sets in polar spaces; many families of strongly regular Cayley graphs with new parameters have been constructed by using cyclotomic classes during the past few years; we survey these constructions and the more recent constructions of $m$-ovoids and $i$-tight sets in classical polar spaces. 


\section{Cyclotomy and difference sets}\label{cycdiff}

The idea of using cyclotomic classes to construct difference sets goes back to Paley~\cite{Paley33}. In the mid-20th century, Baumert, Chowla, Hall, Lehmer, Storer, Whiteman, Yamamoto, etc. pursued this line of research vigorously. Storer's book~\cite{Sto67} contains a summary of results in this direction up to 1967. Important in the study of cyclotomic (or power residue) difference sets are the cyclotomic numbers. Let $q=p^{\ell}$ be a prime power, and let $N>1$ be a divisor of $q-1$. As we did in Section 1, we use $C_i^{(N,q)}, 0\leq i\leq N-1$, to denote the cyclotomic classes of index $N$ of $\F_q$. For integers $a,b$ with $0\leq a,b <N$, the cyclotomic number $(a, b)_N$ 
is defined by
$$(a,b)_N=|(C_a^{(N,q)} +1)\cap C_b^{(N,q)}|.$$

Cyclotomic numbers are useful in many combinatorial investigations, including the study of difference sets in $(\F_q,+)$. These numbers $(a,b)_N$ for $q$ prime have been computed when $N\leq 24$ and $N\not\in\{13, 17, 19, 21, 22, 23\}$ (cf. \cite[p.152]{BEW98}). But it should be noted that when $N$ is large, the formulae given for $(a,b)_N$ are often not explicit. In the following two subsections, we survey recent results on existence/nonexistence results on difference sets in $(\F_q, +)$ arising from unions of cyclotomic classes.

\subsection{A Single Class} We first consider the question when a cyclotomic class $C_i^{(N,q)}$, where $i$ is some integer such that $0\leq i\leq N-1$, is a difference set in $(\F_q,+)$. Since $C_i^{(N,q)}=\gamma^iC_0^{(N,q)}$, the question is equivalent to: When is the cyclotomic class $C_0^{(N,q)}$ a difference set in $(\F_q,+)$? Paley ~\cite{Paley33} is the first to answer this question completely in the case when $N=2$. Later, Chowla~\cite{Cho44} settled the problem in the case when $q$ is prime and $N=4$; Lehmer~\cite{Leh53} gave necessary and sufficient conditions for $C_0^{(N,q)}$ to be a difference set in $(\F_q,+)$ in terms of cyclotomic numbers.

\begin{theorem}\label{lehmer}
Let $C_0^{(N,q)}$ be defined as above. Then $C_0^{(N,q)}$ is a difference set in $(\F_q, +)$ if and only if $N$ is even, $(q-1)/N$ is odd, and 
$$(a, 0)_N=\frac{(q-1-N)}{N^2}$$
for $a=0,1,2,\ldots ,\frac{N}{2} -1$.
\end{theorem}

Theorem~\ref{lehmer} is useful when $N$ is small. Using this theorem, not only one can recover the results of Paley and Chowla, but also obtain complete results in the cases where $N=6$ or $8$.

\begin{theorem}\label{Nsmall}{\em (\cite{Leh53})}
Let $\F_q$ be the finite field of order $q$, where $q=p^{\ell}$ is a power of an odd prime $p$. Let $N\geq 2$ be an even divisor of $q-1$, and $C_0^{(N,q)}$ be the subgroup of $\F_q^*$ of index $N$.  
\begin{itemize}
\item[(1)] When $N=2$, $C_0^{(2,q)}$ is a difference set in $(\F_q,+)$ if and only if $q\equiv 3\pmod 4$.
\item[(2)] When $N=4$, $C_0^{(4,q)}$ is a difference set in $(\F_q,+)$ if and only if $q=p=1+4t^2$ for some odd integer $t$.
\item[(3)] When $N=6$, $C_0^{(6,q)}$ is never a difference set in $(\F_q,+)$.
\item[(4)] When $N=8$, $C_0^{(8,q)}$ is a difference set in $(\F_q,+)$ if and only if $q=p=1+8u^2=9+64v^2$ for some odd integers $u$ and $v$.
\end{itemize}
\end{theorem}

There are a couple of folklore conjectures in this area. It seems difficult to find the exact origin of these conjectures. The third author of the survey was certainly aware of these conjectures many years ago; for example, the stronger conjecture below was mentioned explicitly in~\cite[p.~246]{FX12} and~\cite{Xiang12}. It is quite certain that the history of these conjectures is much longer. The first conjecture is the weaker conjecture.

\begin{conj}\label{weaker}
Let $\F_q$ be the finite field of order $q$, where $q=p^{\ell}$ is an odd prime power. Let $N\geq 2$ be an even divisor of $q-1$, and $C_0^{(N,q)}$ be the subgroup of $\F_q^*$ of index $N$. If $C_0^{(N,q)}$ is a difference set in $(\F_q,+)$, then $N$ must be a power of $2$.  
\end{conj}

The next conjecture is stronger.

\begin{conj}\label{stronger}
Let $\F_q$ be the finite field of order $q$, where $q=p^{\ell}$ is an odd prime power. Let $N\geq 2$ be an even divisor of $q-1$, and $C_0^{(N,q)}$ be the subgroup of $\F_q^*$ of index $N$. If $C_0^{(N,q)}$ is a difference set in $(\F_q,+)$, then $N=2,4$, or $8$.  
\end{conj}

We mention that in a recent paper~\cite{Xia18}, Xia posed essentially the same conjectures as the above folklore conjectures. (It seems that Xia was unaware of the existence of the folklore conjectures above.) Many researchers worked towards settling these conjectures. In the period 1953-1967, the combined work of seven authors showed the nonexistence of difference sets of the form $C_0^{(N,p)}$ in $(\F_p,+)$ for all $8<N<20$, where $p$ is an odd prime; see the book
\cite{Bau71} and \cite[Chapter 5]{BEW98} for references. In 1970, Muskat and Whiteman~\cite{MW70} obtained partial results for the $N=20$ case.  Evans \cite{Evans99} finally finished the $N=20$ case by proving that $C_0^{(20,p)}$ is never a difference set in $(\F_p, +)$, where $p$ is an odd prime. All these nonexistence results were obtained by using Theorem~\ref{lehmer} and cyclotomic numbers. When $N$ is large, Lehmer's theorem is not very useful since the cyclotomic numbers
involved are difficult to compute; instead Gauss sums and Jacobi sums have proved to be more effective. In a recent paper \cite{Xia18}, by using Jacobi sums and extensive Gr{\"o}bner basis computations of certain overdetemined polynomial systems, Xia proved the following theorem.

\begin{theorem}\label{xiathm1}{\em (\cite{Xia18})}
Let $\F_q$ be the finite field of order $q$, where $q=p^{f}$ is an odd prime power. Let $N\geq 2$ be an even divisor of $q-1$, and $C_0^{(N,q)}$ be the subgroup of $\F_q^*$ of index $N$. If $N\leq 22$ and $N\neq 2, 4$ or $8$, then $C_0^{(N,q)}$ is never a difference set in $(\F_q,+)$.  
\end{theorem}

Very recently, Evans and Van Veen~\cite{EV17} proved nonexistence of power residue difference sets in $(\F_p,+)$ for the case where $N = 24$ and $p$ is a prime by computing cyclotomic numbers with the help of a Mathematica program. 

The investigations of the problem when $C_0^{(N,q)}$ is a difference set in $(\F_q,+)$ have also been motivated by questions in finite geometry. A finite projective plane is said to be {\it flag-transitive} if its group of automorphisms acts transitively on the point-line flags. Clearly Desarguesian planes are flag-transitive. Conversely, it is an old and fundamental conjecture in the theory of projective planes, first mentioned in Higman and McLaughlin
\cite{HM61}, that every flag-transitive finite projective plane is Desarguesian. The following theorem, mainly proved by Kantor \cite{Kan87}, relates flag-transitive projective planes to cyclotomic difference sets.

\begin{theorem}\label{flagtrans}
If there exists a non-Desarguesian flag-transitive projective plane of order $n$, then $n^2+n+1:=p$ is prime, $n>8$ is even, and $C_0^{(n,p)}$ is a $(p, n+1, 1)$-difference set in $(\F_p,+)$.
\end{theorem}

By the above theorem, the validity of Conjecture~\ref{stronger} implies that finite flag-transitive projective planes must be Desarguesian. This provided strong motivations to investigate Conjectures~\ref{weaker} and ~\ref{stronger}. Even though many researchers have worked on Conjectures~\ref{weaker} and ~\ref{stronger} for more than sixty years, it seems that we are still far from solving these conjectures. Thas and Zagier
\cite{TZ08} investigated the special case of Conjectures~\ref{weaker} and ~\ref{stronger} related to flag-transitive projective planes. They \cite{TZ08} called a pair $(p,n)$ {\it special}, where $p$ is an odd prime and $1<n<p-1$ an integer dividing $p-1$, if $C_0^{(n,p)}$ is a $(p, n+1, 1)$-difference set in $(\F_p,+)$. Using nontrivial computations, Thas and Zagier \cite{TZ08} classified all special pairs $(p,n)$, when $p<4\times 10^{22}$; no surprises arise from the classification. 

To end this subsection, we caution the readers that two papers with serious mistakes got published during the past 30 years. Feit \cite{Feit90} claimed that if there is a non-Desarguesian projective plane of order $n$, then $n$ is not a power of 2. In \cite{Ott04}, Ott claimed that any flag-transitive finite projective plane has prime power order. Together with Theorem~\ref{flagtrans}, these two results would imply the nonexistence of non-Desarguesian flag-transitive finite
projective planes. Unfortunately both papers, \cite{Feit90} and \cite{Ott04}, contain serious mistakes. We refer the readers to \cite{YH05} and \cite{TZ08} for the exact places in \cite{Feit90, Ott04} where the mistakes were made.

\subsection{Two or More Classes} If Conjecture~\ref{stronger} is true, then $C_0^{(N,q)}$ is rarely a difference set in $(\F_q,+)$. So a natural question is: When is a union of two or more cyclotomic classes a difference set in $(\F_q,+)$ while a single cyclotomic class is not? So far there have been very few results on this question. The first result is a constructive one due to Marshal Hall Jr. \cite{Hall56}. See also \cite[Section 11.6]{Hall86}.

\begin{theorem}
Let $q$ be an odd prime power of the form $q=4x^2+27$ for some integer $x$. Then $C_0^{(6,q)}\cup C_1^{(6,q)}\cup C_3^{(6,q)}$ is a $(q, \frac{q-1}{2}, \frac{q-3}{4})$ difference set in $(\F_q,+)$.
\end{theorem}

The difference sets arising from the above theorem are usually called {\it the Hall sextic residue difference sets}. They were first constructed in the case where $q$ is a prime of the form $4x^2+27$. Later in \cite{Hall86}, the construction was done in the more general setting where $q$ is a prime power of the form $4x^2+27$. However, we note that, as pointed out in \cite{Ocath13}, there are only finitely many proper prime powers of the form $4x^2+27$.  A second remark is
that the above theorem was proved in \cite{Hall56, Hall86} by rather detailed computations of the cyclotomic numbers $(a,b)_6$. It would be interesting to have a proof without using cyclotomic numbers. The reason is that having such a proof will probably pave the way for discovering new difference sets. The investigations of cyclotomic difference sets in the 20th century relied heavily on cyclotomic numbers which are in general very difficult to compute if $N$ is large. It appears that methods using Gauss sums and Jacobi sums directly are more effective for large $N$.

After Marshall Hall Jr.'s work in 1956, several researchers investigated the question when a union of two or more cyclotomic classes is a difference set in the cases where $N=8, 10$, or $12$; only one sporadic difference set, a $(31, 6, 1)$-difference set which is a union of two cyclotomic classes, was found \cite{Hay65} in the case where $N=10$.  Most researchers thought that no new difference sets can be found by taking unions of cyclotomic classes. Therefore it came as a great surprise that in 2012 Feng and the third author \cite{FX12} found new infinite families of difference sets by taking unions of cyclotomic classes with $N=2p_1^m$, where $p_1$ is a prime. We give the detailed statement below. (A difference set $D$ in an additively written finite group $G$ is called {\it skew Hadamard} if $G$ is the disjoint union of $D$, $-D$, and $\{0\}$. A skew Hadamard difference set in a group of order $v$ necessarily has parameter $(v, \frac{v-1}{2}, \frac{v-3}{4})$.)

\begin{theorem}\label{7mod8}{\em (\cite{FX12})}
Let $p_1\equiv 7\, \pmod{8}$ be a prime, $N=2p_1^m$, and 
let $p$ be a prime such that $f:={\rm ord}_N(p)=\phi(N)/2$. Let $s$ be an odd integer, $q=p^{fs}$,  $I$ any subset of $\Z/N\Z$ such that 
$\{i\,\pmod{p_1^m}\,|\,i\in I\}=\Z/{p_1^m}\Z$, and  let $$D=\bigcup_{i\in I}C_i^{(N,q)}\subseteq \F_{q}^*.$$ Then $D$ is a skew Hadamard 
difference set in $(\F_q, +)$  if $p\equiv 3\,\pmod{4}$.
\end{theorem}

Several remarks are in order. First, the proof of the above theorem uses index 2 Gauss sums instead of cyclotomic numbers. Second, the difference sets from Theorem~\ref{7mod8} are not cyclic since the $f$ satisfying the conditions of the theorem is always greater than $1$. Third, there is a lot of flexibility in choosing the index set $I$ in Theorem~\ref{7mod8}; namely, there are $2^{p_1^m}$ choices for the index set $I$ since each pair $\{i, i+p_1^m\}$, $0\leq i\leq p_1^m-1$, contributes exactly one element to $I$. Fourth, the inequivalence between the difference sets from Theorem~\ref{7mod8} and the Paley difference sets was proved by the first author in \cite{Momi13} by using triple intersection numbers.

The case where $p_1$ is a prime congruent to $3$ modulo 8 and $N=2p_1^m$ is more complicated.  Feng and the third author \cite{FX12} first gave a construction of skew Hadamrd difference sets in the case where $N=2p_1$, $p_1\equiv 3\pmod 8$ is a prime. Later on, this construction was generalized by Feng, Momihara and Xiang \cite{FMX15} to work in the case where $N=2p_1^m$, $p_1\equiv 3\pmod 8$ is a prime. Below we state the construction from \cite{FMX15}.

\begin{theorem}\label{3mod8} {\em (\cite{FMX15})}
Let $p_1\equiv 3\pmod{8}$ be a prime, $p_1\neq 3$, $N=2p_1^m$, and 
let $p\equiv 3\pmod{4}$ be a prime such that $f:={\rm ord}_N(p)=\phi(N)/2$. 
Let $q=p^f$, $J=\langle p\rangle \cup 2\langle p\rangle
\cup \{0\}\pmod{2p_1},$
and 
define 
\[
D=\bigcup_{i=0}^{p_1^{m-1}-1}\bigcup_{j\in J}C_{2i+p_1^{m-1}j}
\]  
Assume that $1+p_1=4p^h$, where $h$ is the class number of $\Q(\sqrt{-p_1})$. Then $D$ is a skew Hadamard 
difference set  in the additive group of $\F_q$. 
\end{theorem}

Note that in Theorem~\ref{3mod8}, we need to choose a suitable primitive element $\gamma$ of $\F_q$ in order for the construction to work. We refer the reader to \cite{FMX15} for details on how to choose such a primitive element of $\F_q$.

\section{Sequences with low correlation from cyclotomy}\label{sec-seq}

In this section, we survey results on binary and quaternary sequences with low correlation. Since there exist several excellent surveys on this subject, e.g.~\cite{Ara11,CDR04,GG05,HK98,Schmidt16}, we will concentrate on sequences constructed by using cyclotomy. As indicated in (\ref{eqn-cs}), binary sequences with two-level periodic autocorrelation $\{-1, v\}$ are equivalent to cyclic difference sets  with parameters $(v,(v-1)/2, (v-3)/4)$. Cyclotomy is a powerful tool for constructing
such cyclic difference sets, as we saw in Section~2. Note that the Paley difference set is the classical example of such cyclic difference sets (with $v = p$ a prime) from cyclotomy, and the corresponding characteristic sequence is usually called {\em the Legendre sequence} since the sequence can be defined by the Legendre symbol. In addition, binary sequences of composite length, and quaternary sequences, can also be explicitly constructed using cyclotomy. Below we give a summary of results on such sequences constructed from cyclotomy.

\subsection{Binary sequences from cyclotomy}

By (\ref{eqn-cs}), clearly we have $\mathcal{A}_{\bf s} (\tau) \equiv v \pmod{4}$. Thus, it is natural to classify binary sequences into four categories according to $v \equiv 3 \pmod{4}$, $v \equiv 2 \pmod{4}$, $v \equiv 1 \pmod{4}$, and $v \equiv 0 \pmod{4}$. For each of these four categories, cyclotomy has played an important role in constructing such binary sequences. For $v \equiv 3 \pmod{4}$, binary sequences with two-level autocorrelation $\{-1,
v\}$ are said to have {\em ideal} autocorrelation (for good surveys, see~\cite{CD09,GG05,Xiang99}). It seems very difficult to completely classify binary sequences with ideal autocorrelation, either in terms of sequences or in terms of their supports which are cyclic difference sets. Among the known constructions, there are three arising from cyclotomy:
\begin{itemize}
  \item[(1)] the characteristic sequences of Paley difference sets~\cite{Paley33}; 
  \item[(2)] the characteristic sequences of Hall sextic difference sets~\cite{Hall56};  
  \item[(3)] the twin-prime sequences involving cyclotomic classes of index 2 in both $\F_p$ and $\F_{p+2}$~\cite{SS58}, where $p$ and $p+2$ are twin primes. 
\end{itemize}
We remark that $p$-ary sequences with ideal two-level autocorrelation $\{-1,v\}$ are equivalent to relative difference sets with Singer parameters, and are characterized by the $d$-homogeneous property~\cite{PW151,PW152}. 

A natural question to ask is whether there exist binary sequences with two-level autocorrelation in the other three categories for which $v \not\equiv 3 \pmod{4}$. This question remains open. However, it is evident that the optimal cases for $v \not\equiv 3 \pmod{4}$ are binary sequences with three-level autocorrelation~\cite{JP99} (called {\em optimal} autocorrelation). The supports of such binary sequences with optimal autocorrelation are almost difference
sets. (A subset $D$ of a finite group $G$ is called an {\em almost difference sets} if the list of ``differences'' $d_1 d_2^{-1}$, with $d_1, d_2 \in D$ and $d_1 \ne d_2$ represents each nonidentity element in $G$ either $\lambda$ times or $\lambda + 1$ times~\cite{AD01,DHM01}.) For $v \equiv 2 \pmod{4}$, there are two constructions of binary sequences with three-level autocorrelation $\{2, -2, v\}$ related to cyclotomy: One was given by
Sidelnikov~\cite{Sidel69} (see also~\cite{Turyn68,LCE77}), where the support $D \subseteq (\Z/(q-1)\Z, +)$ is defined as $\log_\gamma (C_1^{(2,q)} - 1)$ with $q \equiv 3 \pmod{4}$ a prime power and $\gamma$ a primitive element of $\F_q$; the other construction was given by Ding, Helleseth and Martinsen~\cite{DHM01}, which in fact uses a union of cyclotomic classes of index $4$ and relies on the explicit computations of cyclotomic numbers. 

For the case $v \equiv 1 \pmod{4}$, all three currently known constructions of binary sequences with autocorrelation values $\{1,-3,v\}$ involve cyclotomy: the first is the Legendre sequence, whose support is the Paley partial difference set; the second was given by Ding, Helleseth and Lam~\cite{DHL99}, and the support is a union of two consecutive cyclotomic classes of index $4$, i.e., $D = C_0^{(4,p)} \cup C_1^{(4,p)}$, where $p = x^2 + 4$ is a prime with $x \equiv 1 \pmod{4}$;
the third construction utilized the so-called generalized cyclotomy, which generalized the twin-prime construction of difference sets to that of almost difference sets by cyclotomic classes of index $2$ in both $\F_p$ and $\F_{p+4}$, where both $p$ and $p+4$ are primes. We note that the second construction $D= C_0^{(4,p)} \cup C_1^{(4,p)}$ was discussed in~\cite{WQWX07}, where the corresponding pseudo-Paley graphs were distinguished from the classical Paley graphs by using $p$-ranks. 

Most of the constructions in the case where $v\equiv 0 \pmod{4}$ interleave four appropriately shifted copies of binary sequences with ideal two-level autocorrelation, while the construction by Sidelnikov~\cite{Sidel69} is an exception: $D:= \log_\gamma (C_1^{(2,q)} - 1)$ with $q \equiv 1 \pmod{4}$ a prime power and $\gamma$ a primitive element in $\F_q$.  

\subsection{Quaternary sequences from cyclotomy}

Given a quaternary sequence ${\bf s}$ of period $v$ over $\{1, i,-1,i^3 \}$ where $i = \sqrt{-1}$, the {\em periodic autocorrelation} at shift $\tau$ with $0 \leq \tau < v$ is defined as ${\mathcal A}_{\bf s} (\tau) = \sum_{i=0}^{v-1} s_i \overline{s_{i+\tau}}$, where $i + \tau$ is read modulo $v$. Each quaternary sequence can be interpreted as two binary sequences via the inverse Gray mapping $\phi^{-1}: \Z/2\Z \times \Z/2\Z \rightarrow \Z/4\Z$, where $\phi^{-1} (0,0) = 0$, $\phi^{-1}
(0,1) = 1$, $\phi^{-1} (1,1) = 2$, and $\phi^{-1} (1,0) = 3$. There are many results on quaternary
sequences with binary sequences with low autocorrelation as building blocks due to~\cite[Eqn. (6)]{KS84}. Instead of giving a complete survey of these results in this section (for recent progress, see for example~\cite{MW18}), we present two constructions of quaternary sequences directly from cyclotomic classes.  

The first construction again is due to Sidelnikov~\cite{Sidel69}, which generates quaternary sequences by $\log_\gamma ( C_j^{(4,q)} - 1)$ for $j = 0,1,2,3$ with $q-1$ divisible by $4$ and $\gamma$ a primitive element in $\F_q$. More generally, for an arbitrary divisor $M$ of $q-1$, $M$-ary sequences of period $q-1$ are obtained in this way with autocorrelation upper bounded by $4$.  

Very recently, a construction of quaternary sequences with autocorrelation bounded by $3$ was proposed in~\cite{MW18} from cyclotomic classes of index $8$. Let $p = x^2 + 16 = a^2 + 2b^2 \equiv 1 \pmod{16}$ ($x \equiv a \equiv 1 \pmod{4}$) be a prime such that $x-a = 4$. Define $D_0 = C_2^{(8,p)} \cup C_6^{(8,p)}$, $D_1 = C_1^{(8,p)} \cup C_3^{(8,p)}$, $D_2 = C_0^{(8,p)} \cup C_4^{(8,p)}$, and $D_3 = C_5^{(8,p)} \cup C_7^{(8,p)}$, and the quaternary sequence ${\bf s}$ of period $p$ is defined by  
$$
s_t = (\sqrt{-1})^j, \qquad \textrm{if $t \in D_j$},
$$
for $j \in \{0,1,2,3\}$ and $s_0 = 1$. Then the quaternary sequence ${\bf s}$ has autocorrelation values $\{-1,-3,3,p\}$. The proof was completed by an explicit computation of cyclotomic numbers of order $8$. Note that the first several primes satisfying the conditions of this construction are $17, 97, 641, 2417, 6577, 14657$.

\section{Strongly regular Cayley graphs from cyclotomy}\label{sec-srg}

A {\it strongly regular graph} srg$(v,k,\lambda,\mu)$ is a simple and undirected graph, neither complete nor edgeless, that has the following properties:

(1) It is a regular graph of order $v$ and valency $k$.

(2) For each pair of adjacent vertices $x,y$, there are $\lambda$ vertices adjacent to both $x$ and $y$.

(3) For each pair of nonadjacent vertices $x,y$, there are $\mu$ vertices adjacent to both $x$ and $y$.\\

Let $\Gamma$ be a (simple, undirected) graph. The adjacency matrix of $\Gamma$ is the $(0,1)$-matrix $A$ with both rows and columns indexed by the vertex set of $\Gamma$, where $A_{xy} = 1$ when there is an edge between $x$ and $y$ in $\Gamma$ and $A_{xy} = 0$ otherwise. A useful way to check whether a graph is strongly regular is by using the eigenvalues of its adjacency matrix. For convenience we call an eigenvalue {\it restricted} if it has an eigenvector which is not a multiple of the all-ones vector ${\bf 1}$. (For a $k$-regular connected graph, the restricted eigenvalues are the eigenvalues different from $k$.)

\begin{theorem}\label{char}
For a simple graph $\Gamma$ of order $v$, neither complete nor edgeless, with adjacency matrix $A$, the following are equivalent:
\begin{enumerate}
\item $\Gamma$ is strongly regular with parameters $(v, k, \lambda, \mu)$ for certain integers $k, \lambda, \mu$,
\item $A^2 =(\lambda-\mu)A+(k-\mu) I+\mu J$ for certain real numbers $k,\lambda, \mu$, where $I, J$ are the identity matrix and the all-ones matrix, respectively, 
\item $A$ has precisely two distinct restricted eigenvalues.
\end{enumerate}
\end{theorem}

For a proof of Theorem~\ref{char}, we refer the reader to \cite{BH12}. An effective method to construct strongly regular graphs is by using Cayley graphs. Let $G$ be an additively written group of order $v$, and let $D$ be a subset of $G$ such that $0\not\in D$ and $-D=D$, where $-D=\{-d\mid d\in D\}$. The {\it Cayley graph on $G$ with connection set $D$}, denoted by ${\rm Cay}(G,D)$, is the graph with the elements of $G$ as vertices; two vertices are adjacent if and only if their difference belongs to $D$. In the case when $\Cay(G,D)$ is a strongly regular graph, the connection set $D$ is called a (regular) {\it partial difference set}. Examples of strongly regular Cayley graphs are the Paley graphs ${\rm P}(q)$, where $q$ is a prime power congruent to 1 modulo 4, the Clebsch graph, and the affine orthogonal graphs (\cite{BH12}). For $\Gamma={\rm Cay}(G,D)$ with $G$ abelian, the eigenvalues of $\Gamma$ are exactly $\chi(D):=\sum_{d\in D}\chi(d)$, where $\chi$ runs through the character group of  $G$. This fact reduces the problem of computing eigenvalues of abelian Cayley graphs to that of computing some character sums, and is the underlying reason why the Cayley graph construction has been very effective for the purpose of constructing srgs. The survey of Ma~\cite{Ma94} contains much of what is known about partial difference sets and about connections with strongly regular graphs. 

In this section, we always take the additive group of a finite field as the underlying group $G$ and take a union of cyclotomic classes as connection sets.  Many reseachers have studied the problem of determining when a union $D$ of cyclotomic classes forms a partial difference set.  In some of the papers, the authors used the language of codes or finite geometry in their studies instead of strongly regular Cayley graphs or partial difference sets. We choose to use the language of srgs here.

\begin{example}\label{ex:spor}{\em 
Here are three known ``sporadic'' examples of strongly regular Cayley graphs on finite fields: 
\begin{itemize}
  \item[(1)] (\cite{vS81}) $\Cay(\F_{3^4},D)$ with $D=\bigcup_{i\in \{0,1,3\}}C_i^{(8,3^4)}$ is an srg$(3^4,30,9,12)$; 
\item[(2)] (\cite{Hill76}) $\Cay(\F_{2^{12}},D)$ with $D=\bigcup_{i\in \{0,7\}}C_i^{(35,2^{12})}$ is an srg$(2^{12},234,2,14)$; 
\item[(3)] (\cite{de95}) $\Cay(\F_{3^8},D)$ with $D=\bigcup_{i\in \{0,1,2,8,10,11,13\}}C_i^{(16,3^8)} $ is an srg$(3^8, 2870,1249,1260)$. 
\end{itemize}
}
\end{example}

\subsection{Cyclotomic strongly regular graphs}\label{sec:CSRG}
Let $p$ be a prime,  $\ell$ and $m$ be positive integers, and let $q=p^\ell$. Let $N>1$ be an integer such that $N|(q^m-1)$, and $\gamma$ be a primitive element of $\F_{q^m}$. For a subset $D$ of $\F_{q^m}^*$, 
we call $\Cay(\F_{q^m},D)$ a {\it cyclotomic strongly regular graph} if $D$ is a single cyclotomic class of $\F_{q^m}$ and $\Cay(\F_{q^m},D)$ is strongly regular. 
The Paley graphs are primary examples of  cyclotomic srgs.  Also, if $D$ is the multiplicative group of a subfield of $\F_{q^m}$, then it is clear that $\Cay(\F_{q^m} , D)$ is strongly regular.  These cyclotomic srgs are usually called {\it subfield examples}. Next, if there exists a positive integer $j$ such that $p^j\equiv -1\,(\mod{N})$, then $\Cay(\F_{q^m}, D)$ is strongly regular. See \cite{BMW82} for a proof of this result. These examples are usually called {\it semi-primitive}.
A generalization of semi-primitive srgs so that its connection set is a union of at least two cyclotomic classes 
was given in \cite{BWX99}; that generalization will be explained in Subsection~\ref{subsec:semi}.

In \cite{SW02}, Schmidt and White gave the following necessary and sufficient condition for $\Cay(\F_{q^m},D)$ to be a cyclotomic srg. 
\begin{theorem}\label{SW_NScond}{\em (\cite{SW02})}
With notation as above, assume that  $N$ divides $(q^m-1)/(q-1)$.  Let $f$ be the order of $p$ modulo $N$, and put $s=m\ell/f$. 
Then, $\Cay(\F_{q^m},C_0^{(N,q^m)}) $ is strongly regular if and 
only if there exists a positive integer $u$ satisfying the following three conditions: 
\begin{itemize}
\item[(i)]  $u\,|\,(N-1)$; 
\item[(ii)] $up^{st}\equiv \pm 1\,(\mod{N})$;
\item[(iii)] $u(N-u)=(N-1)p^{s(f-2t)}$. 
\end{itemize}
Here, $t$ is the largest power of $p$ dividing the Gauss sums  $G_{q^m}(\chi)$ for all nontrivial multiplicative character $\chi$ of $\F_{q^m}$ of order dividing $N$. 
\end{theorem}
The necessary and sufficient conditions in the above theorem can be used to search for cyclotomic srgs $\Cay(\F_{q^m}, C_0^{(N, q^m)})$ with large $N$. The eleven sporadic examples in Table~\ref{Tab1} which are neither subfield examples nor semi-primitive examples were found in this way in~\cite{SW02} (some of the eleven examples in Table~\ref{Tab1} were already known before the search conducted in \cite{SW02}; see \cite{BaD99, Lang96}). A generalization of these sporadic examples so that their connection sets are union of at least two cyclotomic classes 
was given in \cite{FMX15,FX122,GXY13,Momi132}.  We will explain that generalization in Subsection~\ref{subsec:spo}.

\begin{table}[h]
\caption{Eleven sporadic examples}
\label{Tab1}
$$
\begin{array}{|c||c|c|c|c|}
\hline
\mbox{No.}&N&q&m&[(\Z/N\Z)^\ast:\langle p\rangle]\\
\hline
1&11&3&5&2\\
2&19&5&9&2\\
3&35&3&12&2\\
4&37&7&9&4\\
5&43&11&7&6\\
6&67&17&33&2\\
7&107&3&53&2\\
8&133&5&18&6\\
9&163&41&81&2\\
10&323&3&144&2\\
11&499&5&249&2\\
\hline
\end{array}
$$
\end{table}

On the other hand, Schmidt and White~\cite{SW02} made the following conjecture on cyclotomic srgs, which can be thought as a counterpart of Conjecture~\ref{stronger} for cyclotomic srgs. 

\begin{conj}\label{con:SW}{\em (\cite{SW02})}
Assume that $N\,|\,(q^m-1)/(q-1)$. 
Then, $\Cay(\F_{q^m},C_0^{(N,q^m)}) $ is strongly regular if and 
only if  it is either a subfield example, or a semi-primitive example or one of the eleven sporadic 
examples of Table~\ref{Tab1}. 
\end{conj}
The Schmidt-White conjecture remains open. There are some results on this conjecture in~\cite{SW02} under the condition $[(\Z/N\Z)^\ast:\langle p\rangle]=2$ and the assumption of the generalized Riemann hypothesis. 

\begin{remark}
Theorem~\ref{SW_NScond} and Conjecture~\ref{con:SW} were stated  in terms of  two-weight  irreducible cyclic codes in  \cite{SW02}. 
We briefly explain the connection between two-weight  irreducible cyclic codes and cyclotomic srgs below.  

For a positive divisor $n$ of $q^m-1$, let $\xi$ be a primitive $n$th root of unity in $\F_{q^m}$. Then,   
\[
C=\left\{c(y):=\left(\Tr_{q^m/q}(y \xi^i)\right)_{i=0}^{n-1}\,|\,y\in \F_{q^m}\right\} 
\]
is called an {\it irreducible cyclic code} of length $n$ over $\F_q$. McEliece~\cite{Mc74} showed that if $N:=(q^m-1)/n$ divides $(q^m-1)/(q-1)$, 
the Hamming weight of $c(y)$ for $y\in \F_{q^m}^\ast$ is given by 
\[
\frac{(q-1)}{qN}(q^m-1-N\cdot \psi_{\F_{q^m}}(yC_0^{(N,q^m)})), 
\] 
where $\psi_{\F_{q^m}}$ is the canonical additive character of $\F_{q^m}$. 
Hence, $C$ is a two-weight code if and only if 
$\psi_{\F_{q^m}}(yC_0^{(N,q^m)})$, $y\in \F_{q^m}^\ast$, take exactly two values, i.e., $\Cay(\F_{q^m},C_0^{(N,q^m)})$ is strongly regular. 
For more details on the correspondence between projective two-weight codes and strongly regular Cayley graphs on finite fields, see, e.g., \cite[p.~140]{BH12}.  
\end{remark}

\subsection{A generalization of semi-primitive examples}\label{subsec:semi}
Let $q=p^m$ be a prime power with $p$ a prime and $N$ be a 
positive integer dividing $q-1$. Let $\gamma$ be a primitive element of $\F_q$. 
 Assume that 
there is a $j>0$ such that $p^j\equiv -1\,(\mod{N})$. 
Choose $j$ minimal with this property and write $m=2js$. 

The following theorem is a generalization of semi-primitive examples of cyclotomic srgs so that their connection sets are unions of at least two cyclotomic classes.  
\begin{theorem}\label{thm:semip}{\em (\cite{BWX99,CK86})}
With notation as above, let $J$ be a subset of $\{0,1,\ldots,N-1\}$ of size $\ell$ and 
$D=\bigcup_{i\in J}C_i^{(N,q)}$. If $D=-D$, then 
$\Cay(\F_{q},D)$ is an srg with parameters $(u^2,r(u-\epsilon)\epsilon u+r^2-3\epsilon r,r^2-\epsilon r)$ with 
$u= p^{js}$ and $r=\ell(p^{js}+\epsilon)/N$, where $\epsilon=-1$ or $1$  
depending on whether $s$ is even or odd. 
In particular, for $a=0,1,\ldots,N-1$,  
\[
\psi_{\F_{q}}(\gamma^a D)=\frac{u((-1)^s \sqrt{q}-1)}{N}+ 
\begin{cases}
(-1)^{s+1}\sqrt{q},& \text{ if $\delta^s=1$ and $a\in -J\,(\mod{N})$} \\
& \text{ \, \, or $\delta^s=-1$ and $a\in -J+N/2\,(\mod{N})$}, \\
0, & \text{ otherwise, }
\end{cases}
\]
where \[
\delta=\begin{cases}
1,& \text{ if $N$ is even and $(p^j+1)/N$ is odd}, \\
-1, & \text{ otherwise.}
\end{cases}
\]
\end{theorem}

We mention that an srg is said to be of {\it Latin square type} (respectively, {\it negative Latin square type}) if $(v,k,\lambda,\mu)=
(u^2,r(u-\epsilon),\epsilon u+r^2-3\epsilon r,r^2-\epsilon r)$ and $\epsilon=1$ (respectively, $\epsilon=-1$). Most known strongly regular Cayley graphs are of Latin square or negative Latin square type.

In~\cite{BWX99}, the following two further generalizations were given. Pick several positive integers $N_i$, $i\in I$, with $N_i\,|\,(q-1)$.  For each $i \in I$, let $J_i$ be a subset of $\{0,1,\ldots,N_i-1\}$. We define $D_{i}=\bigcup_{j \in J_i}C_{j}^{(N_i,q)}$, and assume that $D_i$ are mutually disjoint. Then, it is possible for $D$ to give rise to a strongly regular Cayley graph. The precise statements are given below.

\begin{proposition}\label{prop:ex1}
Let $p$ be an odd prime and $N_i=p^{j_i}+1$ for $i=1,2$. Let $q=p^m$ with $m=4j_1 s_1=4j_2s_2$. Take $J_1$ as a subset of $\{2h\,|\,h=0,1,\ldots,N_1/2-1\}$ and  $J_2$ as a subset of $\{2h+1\,|\,h=0,1,\ldots,N_2/2-1\}$. Define $D_{i}=\bigcup_{j \in J_i}C_{j}^{(N_i,q)}$, $i=1,2$, and $D=D_1\cup D_2$. Then, $\Cay(\F_{q},D)$ is an srg of negative Latin square type.   
\end{proposition}

The proof of  Proposition~\ref{prop:ex1}  is obvious since $D_1\cap D_2=\emptyset$ and $a\in -J_i\,(\mod{N_i})$ cannot hold for $i=1$ and $2$ simultaneously. 

\begin{example}\label{exm:DeLange1}{\em 
Let $(p,j_1,j_2,N_1,N_2,f)=(3,1,2,4,10,8)$, $J_1=\{0\}$, $J_2=\{1\}$, and $J=\{0,1,4,8,11,12,16\}$. Then, 
\[
D=\bigcup_{i=1,2}\bigcup_{h\in J_1}C_h^{(N_i,q)}=\bigcup_{h\in J}C_h^{(20,q)},  
\]
and $\Cay(\F_{q},D)$ is an srg with parameters 
$(v,k,\lambda,\mu)=(3^8,2296,787,812)$. }
\end{example}

Similar to Proposition~\ref{prop:ex1}, we have the following. 

\begin{proposition}\label{prop:ex2}
Let $p$ be an odd prime and $N_i=p^{j_i}+1$ for $i=1,2$. Let $q=p^m$ with $m=2j_1 s_1=2j_2s_2$, where $s_1$ and $s_2$ are odd. Take $J_1$ as a subset of $\{2h\,|\,h=0,1,\ldots,N_1/2-1\}$ and  $J_2$ as a subset of $\{2h+1\,|\,h=0,1,\ldots,N_2/2-1\}$. Assume that $\bigcup_{i\in J_1}C_{-i+N_1/2}^{(N_1,q)}$ and $\bigcup_{i\in J_2}C_{-i+N_2/2}^{(N_2,q)}$ are disjoint. Define $D_{i}=\bigcup_{j \in J_i}C_{j}^{(N_i,q)}$, $i=1,2$, and $D=D_1\cup D_2$. Then, $\Cay(\F_{q},D)$ is an srg of negative Latin square type.   
\end{proposition}

There are many choices of $p$ and $j_i$, $i=1,2$, 
satisfying the condition of Proposition~\ref{prop:ex2}. For example, 
if $p\equiv 3\,(\mod{4})$ and $j_1$ and $j_2$ are both odd, then 
$\bigcup_{i\in J_1}C_{-i+N_1/2}^{(N_1,q)}$ and $\bigcup_{i\in J_2}C_{-i+N_2/2}^{(N_2,q)}$ are disjoint. 
\subsection{A generalization of sporadic or subfield examples}\label{subsec:spo}

In~\cite{FMX15,FX122,GXY13}, the authors found infinite families of strongly regular Cayley graphs on finite fields generalizing seven of the eleven sporadic examples of cyclotomic srgs in Table~\ref{Tab1}. Their constructions used unions of ``consecutive'' cyclotomic classes of finite fields as connection sets for the Cayley graph construction. In particular, the following theorem was proved.

\begin{theorem} \label{thm:knows}
\begin{itemize}
\item[(i)] {\em (\cite{FX122})} Let $q=p^{p_1^{m-1}(p_1-1)/2}$, $N=p_1^m$, and  $D=\bigcup_{i=0}^{p_1^{m-1}-1}C_{i}^{(N,q)}$. Then, $\Cay(\F_q,D)$ is strongly for any $m\ge 1$ in 
the following cases:  
\begin{eqnarray*}
(p,p_1)=(2,7),(3,107),(5,19),(5,499),(17,67),(41,163). 
\end{eqnarray*}
\item[(ii)] {\em (\cite{GXY13})} Let $q=p^{p_1^{m-1}(p_1-1)/4}$, $N=p_1^m$, and  $D=\bigcup_{i=0}^{p_1^{m-1}-1}C_{i}^{(N,q)}$. 
Then, $\Cay(\F_q,D)$ is strongly regular for any $m\ge 1$ in 
the following cases:  
\begin{eqnarray*}
(p,p_1)=(3,13),(7,37). 
\end{eqnarray*}
\item[(iii)] {\em (\cite{FMX15})} Let $q=p^{p_1^{m-1}(p_1-1)p_2^{n-1}(p_2-1)/2}$, $N=p_1^mp_2^n$, and  
$D=\bigcup_{i=0}^{p_1^{m-1}-1}\bigcup_{j=0}^{p_2^{n-1}-1}
C_{p_2^n i+p_1^m j}^{(N,q)}$. Then, $\Cay(\F_q,D)$ is strongly regular for any $m,n\ge 1$ in 
the following cases:  
\begin{eqnarray*}
(p,p_1,p_2)=(2,3,5),(3,5,7),(3,17,19). 
\end{eqnarray*}
\end{itemize}
\end{theorem}

The srgs in the cases when $(p,p_1)=(2,7),(3,13)$ and $(p,p_1,p_2) = (2,3,5)$ in  Theorem~\ref{thm:knows} are generalizations of subfield examples. The others are generalizations of sporadic examples of Table~\ref{Tab1}. In all cases, 
it holds that $[(\Z/N\Z)^\ast:\langle p\rangle]=[(\Z/p_1\Z)^\ast:\langle p\rangle]$ or $[(\Z/N\Z)^\ast:\langle p\rangle]=[(\Z/p_1p_2\Z)^\ast:\langle p\rangle]$. The proofs are based on known evaluations of index $2$ or $4$ Gauss sums (see~\cite{FYL05,YX10}). 

Note that it is unlikely that one can generalize the 1st example in Table~\ref{Tab1} by a similar method since $[(\Z/11^m\Z)^\ast:\langle 3\rangle]\not=[(\Z/11\Z)^\ast:\langle 3\rangle]$ for $m\ge 2$. 
In order to generalize the 5th and 8th srgs in Table~\ref{Tab1} into infinite families, we may need to evaluate Gauss sums of index $6$. However, it seems very difficult to compute Gauss sums of index  $e$ when $e>4$. As a result, it is hard to find new srgs on $\F_q$ in the index $e>4$ cases. On the other hand, in \cite{Momi132}, the first author of this survey succeeded in giving a recursive construction of srgs, which enables him to generalize the remaining examples into infinite families not using explicit evaluations of Gauss sums. Instead, he studied the rationality of  ``relative'' Gauss sums. 

\begin{theorem}\label{thm:srg1}{\em (\cite{Momi132})}
Let $N_1=p_1\cdots p_mp_{m+1}\cdots p_\ell$, where $p_i$'s are distinct odd primes, and assume that $[(\Z/h\Z)^\ast:\langle p\rangle]=e$.  Furthermore, Let $N=p_1^{e_1}\cdots p_m^{e_m} p_{m+1}^{e_{m+1}}\cdots p_\ell^{e_{\ell}}$, where $e_i\ge 1$ for $1\le i\le m$ and $e_i= 1$ for $m+1\le i\le \ell$, and assume that $\langle p\rangle $ is of index $e$ modulo $N$. Let $q_1=p^{d}$ and $q=p^{f}$, where $d=\phi(N_1)/e$ and  $f=\phi(N)/e$. Here, $\phi$ is the Euler totient function. Put $h_j=\prod_{i\not=j} p_i$ for $1\le j\le m$. Assume that there exists an integer $s_j$ such that $p^{s_j}\equiv -1\,(\mod{h_j})$ for $1\le j\le m$. Let 
\[
D:=\bigcup_{i_1=0}^{p_1^{e_1-1}-1}\cdots \bigcup_{i_m=0}^{p_m^{e_m-1}-1}C_{i_1{n}_1+\cdots+i_m {n}_m }^{(N,q)}, 
\]  
where $n_j=\prod_{i\not=j} p_i^{e_i}$. If $\Cay(\F_{q_1},C_0^{(N_1,q_1)})$ is an srg, then  so is $\Cay(\F_q,D)$.
\end{theorem}

\begin{example}{\em 
\begin{itemize}
\item[(i)] We can apply Theorem~\ref{thm:srg1} to the 5th srg in Table~\ref{Tab1} as $(\ell,p_1,p,e)=(1,43,11,6)$. In this case,  we do not need the condition that there exists an integer $s_j$ such that $p^{s_j}\equiv -1\,(\mod{h_j})$. It is clear that $[(\Z/p_1^{e_1} \Z)^\ast:\langle p\rangle]=6$ for any $e_1\ge 1$. Hence, $\Cay(\F_{p^{p_1^{e_1-1}(p_1-1)/6}},D)$ is strongly regular, where 
\[
D=\bigcup_{i=0}^{p_1^{e_1-1}-1}C_i^{(p_1^{e_1},p^{p_1^{e_1-1}(p_1-1)/6})}.
\] 
There are many examples in the subfield case satisfying the condition of Theorem~\ref{thm:srg1} with $\ell=1$, for example,    
\[
(p,f,p_1,e)=(3,3,13,4), (2,5,31,6),(5,3,31,10),(2,9,73,8). 
\] 
In these cases, we have $[(\Z/p_1^{e_1} \Z)^\ast:\langle p\rangle]=e$ for any $e_1\ge 1$. Hence, these examples can be similarly generalized into infinite families. 
\item[(ii)] 
We can apply
Theorem~\ref{thm:srg1} to the 8th srg in Table~\ref{Tab1} as $(\ell,m,p_1,p_2,p,e)=(2,1,19,7,5,6)$. In this case, there exists an integer $s_2$ such that $p^{s_2}\equiv -1\,(\mod{p_2})$. It is clear that $[(\Z/p_1^{e_1}p_2 \Z)^\ast:\langle p\rangle]=6$ for any $e_1\ge 1$. Hence, \linebreak
$\Cay(\F_{p^{p_1^{e_1-1}(p_1-1)(p_2-1)/6}},D)$ is strongly regular, where 
\[
D=\bigcup_{i=0}^{p_1^{e_1-1}-1}C_i^{(p_1^{e_1}p_2,p^{p_1^{e_1-1}(p_1-1)(p_2-1)/6})}.
\] 
There are many examples in the subfield case satisfying the condition of Theorem~\ref{thm:srg1} with $\ell=2$, for example,    
\[
(p,f,p_1,p_2,e)=(2,4,3,5,2), (2,8,5,17,8),(2,10,31,11,30),(2,14,127,43,378). 
\] 
In the former two cases, we have 
$[(\Z/p_1^{e_1}p_2^{e_2} \Z)^\ast:\langle p\rangle]=e$ for any $e_1,e_2\ge 1$ 
and $p$ is semi-primitive modulo both $p_1$ and $p_2$. In the latter  
two cases, we have 
$[(\Z/p_1^{e_1}p_2 \Z)^\ast:\langle p\rangle]=e$ for any $e_1\ge 1$, 
and $p$ is semi-primitive modulo $p_2$ only. 
Hence, 
these examples can be generalized into infinite families by using Theorem~\ref{thm:srg1}. 
\end{itemize}}
\end{example}

\subsection{On de Lange's sporadic examples of srgs}\label{subsec:DeLange} 
In \cite{de95}, de Lange found four ``sporadic'' examples of strongly regular Cayley graphs on the additive groups of finite fields by using a computer. The srgs he found have the following parameters: 
\begin{itemize}
\item[(1)] $(v,k,\lambda,\mu)=(3^8,2296,787,812)$; 
\item[(2)] $(v,k,\lambda,\mu)=(3^8,2870,1249,1260)$; 
\item[(3)] $(v,k,\lambda,\mu)=(2^{12},273,20,18)$; 
\item[(4)] $(v,k,\lambda,\mu)=(2^{12},1911,950,840)$.  
  \end{itemize}

The 3rd and 4th examples of srgs are dual to each other; hence de Lange found essentially three examples. In particular, the 2nd example is the one given in Example~\ref{ex:spor}~(3). As explained in Example~\ref{exm:DeLange1} and Theorem~\ref{thm:knows}~(iii), 
the 1st and 3rd examples above have already been generalized in \cite{BWX99} and \cite{FX122}, respectively. However, it seems difficult to generalize the 2nd example above into an infinite family of srgs. In \cite{Xiang12}, the third author asked the question of generalizing the last example of de Lange (see Problem 5.2 in \cite{Xiang12}). In this subsection, we show that there is an infinite family of srgs including an srg with the same parameters as those of the 2nd example above.

We will need the following families of srgs.


\begin{theorem}{\em (\cite{CK86})}\label{thm:affine}
Let $Q:\F_q^n\to \F_q$ be a nonsingular quadratic form, where $n=2m$ is even and $q$ is an odd prime power. Define 
$\cQ=\{x\in \F_q^n\setminus \{{\bf 0}\}\,|\,
Q(x)=0\}$
and 
$D_i=\{x\in \F_q^n\,|\,
Q(x)\in C_i^{(2,q)}\}$, $i=0,1$. 
Then, each $\Cay(\F_q^n,D_i)$, $i=0,1$, is 
an srg with parameters $(u^2,r(u-\epsilon),\epsilon u+r^2-3\epsilon r,r^2-\epsilon r)$ with 
$u = q^{m}$ and $r=\epsilon q^{m-1}(q-1)/2$, where $\epsilon=1$ or $-1$  
depending on whether $Q$ is hyperbolic or  elliptic. 
\end{theorem}  

The srg $\Cay(\F_{q}, D_0)$ in the above theorem is called an {\it affine polar graph}. In \cite{MX14}, the following recursive construction of srgs was given as a generalization of the construction above. Let $q$ be a prime power and $N>1$ be an integer dividing $q-1$. Furthermore, let $\gamma$ be a fixed primitive element of $\F_{q^2}$, and let  $\omega=\Norm_{q^2/q}(\gamma)$, which is a primitive element of $\F_q$. Put $C_{i}^{(e,q^2)}=\gamma^i\langle
\gamma^{N}\rangle$, $i=0,1,\ldots,N-1$. Let $Q:\F_q^n \to \F_q$ be a quadratic form. For $y \in \F_q$, define $D_y = \{x\in \F_q^n\,|\,Q(x)=y\}$, and for a subset $E$ of $\F_q$, write  $D_E=\sum_{y \in E}D_y$. 

\begin{theorem}\label{thm:liftQuad}  
Let $J$ be a subset of $\{0,1,\ldots,N-1\}$, and let   $E=\bigcup_{i\in J}C_i^{(N,q)}$. 
Assume that $\Cay(\F_{q^2},\bigcup_{i\in J}C_i^{(N,q^2)})$ is an srg of negative Latin square type. 
Let $Q:\F_q^n \to \F_q$ be a nonsingular quadratic form, where $n=2m$ is even.
Then, $\Cay(\F_q^n,D_{E})$  is an srg$(u^2,r(u-\epsilon),\epsilon u+r^2-3\epsilon r,r^2-\epsilon r)$ with 
$u = q^m$ and $r=\epsilon |J|  q^{m-1}(q-1)/N$, where $\epsilon=1$ or $-1$  
depending on whether $Q$ is hyperbolic or  elliptic. 
%
\end{theorem}

In \cite{MX14,MX18}, the authors gave constructions of strongly regular graph $\Cay(\F_{q^2},D)$ of negative Latin square type such that 
$D$ is a union cosets of $C_0^{(q-1,q^2)}$ based on cyclotomic srgs, which can be used as starters 
in order to apply Theorem~\ref{thm:liftQuad}. Furthermore, in \cite{FWXY13}, the authors studied a construction of srgs based on weakly regular bent functions instead of quadratic forms in Theorem~\ref{thm:liftQuad}.

We will also need the following family of srgs.  

\begin{example}\label{exm:1movoid}{\em 
Let $q$ be an odd prime power and $Q:\F_q^8\to \F_q$ be an elliptic quadratic form defined by 
\begin{equation}\label{eq:ell}
Q(x)=\Tr_{q^8/q}(x^{q^4+1}). 
\end{equation}
We define $\cQ=\{x\in \F_q^n\setminus \{{\bf 0}\}\,|\,
Q(x)=0\}$.
Let $D=C_{(q^2+1)/2}^{(q^2+1,q^8)}$. Then, $D\subseteq \cQ$, and 
$\Cay(\F_{q^8},D)$ is strongly regular with negative Latin square type parameters $(q^8,r(q^4+1),-q^4+r^2+3 r,r^2+r)$, where $r=q^2-1$, since $q^2\equiv -1\pmod N$ with $N=q^2+1$ (i.e., the semi-primitive condition holds here). In this case, $\Cay(\F_{q^8},\cQ\setminus D)$ is also an srg  of the same type.
}\end{example}

Finally, we use the following theorem of van Dam~\cite{vD03}.  
\begin{theorem}\label{vD}
Let $\{G_1,G_2,\ldots,G_d\}$ be a decomposition of the complete graph on a vertex set $X$, where each $G_i$ is strongly regular. If the $G_i$'s are all of Latin square type or all of negative Latin square type, then a union of any subset of $\{G_1,G_2,\ldots,G_d\}$ is also an srg of the same type on $X$. 
\end{theorem}
\begin{remark}
In \cite{vD03}, van Dam actually proved that the decomposition $\{G_1,G_2,\ldots,G_d\}$ forms a $d$-class amorphic association scheme under the same assumption of Theorem~\ref{vD}. We will not need the full strength of this result. The theorem above suffices for our purpose. 
\end{remark}

\begin{example}\label{exm:DeLange3}
{\em Let $q$ be an odd prime power and $n=8$. 
Let $Q:\F_q^n\to \F_q$ be a nonsingular elliptic quadratic form.  Let 
$G_1=\Cay(\F_{q}^n,D_1)$ be an 
affine polar graph of negative Latin square type associated with $Q$, and 
$G_2=\Cay(\F_{q}^n,D_2)$ be the srg defined in \eqref{eq:ell}. Then, the four graphs 
$G_1$, $G_2$, $G_3=\Cay(\F_{q}^n,\F_q^n\setminus (\{{\bf 0}\}\cup \cQ\cup D_1))$, 
$G_4=\Cay(\F_{q}^n,\cQ\setminus D_2)$ give a decomposition of the complete graph on $\F_{q}^n$, where each $G_i$ 
is an srg of negative Latin square type. By Theorem~\ref{vD}, the graph 
$\Gamma=G_1+G_3$ is also an srg of negative Latin square type. 
Take $q=3$, and then the graph $\Gamma$ is an srg with 
parameters $(v,k,\lambda,\mu)=(3^8,2870,1249,1260)$. This srg has the same parameters as those of de Lange's 2nd example $\Gamma'$ of srgs. We checked that $\Gamma$ and $\Gamma'$ are nonisomorphic by a computer. In particular, we have 
$\#\Aut(\Gamma)=2^7\cdot 3^{12}\cdot 5\cdot 41$ and $\#\Aut(\Gamma')=2^4\cdot 3^8\cdot 5\cdot 41$.}  
\end{example}

\subsection{Projective two-intersection sets, $m$-ovoids, and $i$-tight sets}\label{subsec:movoid}
During the past few years, strongly regular Cayley graphs defined on the additive groups of finite fields have been extensively studied due to their close connections with certain substructures in finite geometry. In most published works by geometers, the authors used the language of projective two-intersection sets, or two-character sets.  Because of the large amount of papers published in this direction, it is difficult to summarize all known constructions and existence results in
this short subsection.  Instead we will focus on explaining the connections between projective two-intersection sets and strongly regular Cayley graphs on finite fields, and a linkage with geometric objects, called $m$-ovoids and $i$-tight sets in polar spaces.  

A set ${\mathcal M}$ of points of a projective space $\PG(n-1,q)$  is called a 
{\it projective two-intersection set of type $(a,b)$} (or simply, a {\it set of type $(a,b)$}) if every hyperplane of $\PG(n-1,q)$ meets ${\mathcal M}$ in $a$ or $b$ points. In some papers, a projective two-intersection set is also called a {\it two-character set}.

\begin{example}\label{exam}{\em 
\item[(1)] A hyperoval in $\PG(2,2^f)$ is a set of type  $(0,2)$. 
\item[(2)] A unital in $\PG(2,q^2)$ is a set of type $(1,q+1)$.  
\item[(3)] A nodegenerate quadric $\cQ$ in $\PG(2m-1,q)$ is a set of type $(\theta_{m}-q^{2m-2},\theta_m-q^{2m-2}-\epsilon q^{m-1})$, where $\theta_m=\frac{(q^m-\epsilon)(q^{m-1}+\epsilon)}{q-1}$ and $\epsilon =1$ or 
$-1$ depending on whether  $\cQ$ is  hyperbolic or elliptic. }
\end{example}
See \cite{Den69} for a generalization of (1) in Example~\ref{exam} and \cite{Bro85} for a difference of two quadrics construction.  

Let $N=(q^{n}-1)/(q-1)$, and let $\gamma$ be a primitive element of $\F_{q^{n}}$. We identify the points of $\PG(n-1,q)$ with $\Z_N$ as follows: 
View $\F_{q^{n}}$ as an $n$-dimensional space over $\F_q$, and use 
$\F_{q^{n}}$ as the underlying vector space of $\PG(n-1,q)$. We identify the projective point $\langle \gamma^i\rangle$ with $i\in \Z_N$. 
Then, all hyperplanes in $\PG(n-1,q)$ are given by 
\[
H_i:=\{\langle \gamma^j\rangle\,|\,\Tr_{q^{n}/q}(\gamma^{i+j})=0,j \in \Z_N\}, \, \, 
i \in \Z_N. 
\]
Now let $\mathcal{M}$ be a set of points of $\PG(n-1,q)$,  and define 
\[
D:=\{xy: y \in \F_q^\ast,\langle x\rangle \in {\mathcal M}\}\subseteq \F_{q^n}. 
\]
Then, we have 
\[
\psi_{\F_{q^n}}(\gamma^i D)=
\sum_{y\in \F_{q}}\sum_{x \in {\mathcal M}}\zeta_p^{\tr_{q^n/q}(\gamma^i xy)}-|{\mathcal M}|
=
q |H_i\cap {\mathcal M}|-|{\mathcal M}|. 
\]
Hence, ${\mathcal M}$ is a set of type $(a,b)$ in $\PG(n-1,q)$ if and 
only if the character values of $D$ take exactly two values $qa-|{\mathcal M}|$ and 
$qb-|{\mathcal M}|$, i.e., 
$\Cay(\F_{q^{n}},D)$ is strongly regular with parameters $(q^n, (q-1)|{\mathcal M}|, \lambda, \mu)$, where $\lambda$ and $\mu$ can be computed from $a$, $b$, $|{\mathcal M}|$, $q$, and $n$. 

There are many known constructions of projective two-intersection sets. 
See, e.g., \cite{Cos10,CDMPS08,CK10,CP15,CP13,DV08,De16,IZZ15,Nap15,Nap13,Pav15,PR95},  for recent constructions of projective two-intersection sets. 

Many projective two-intersection sets arise from $m$-ovoids and $i$-tight sets in classical polar spaces. Conversely, projective two-intersection sets with certain special properties can give rise to $m$-ovoids and $i$-tight sets. Many recent constructions of $m$-ovoids and $i$-tight sets came about via constructions of projective 2-intersection sets satisfying special properties, see, e.g. \cite{FMX152, BLMX18}.

Let $V=\F_q^n$ be an $n$-dimensional vector space over $\F_q$ and 
$f$ be a non-degenerate sesquilinear or non-singular quadratic form defined on $V$. A finite classical polar space associated with the form 
$f$ is the geometry consisting of subspaces of $\PG(n-1,q)$ induced by 
the totally isotropic subspaces with relation to $f$. A polar space $S$ contains 
totally isotropic points, lines, planes, etc. The (totally isotropic) subspaces of maximum dimension are called {\it maximals} of $S$. The {\it rank} of $S$ is the vector dimension  of its maximals. 

There are three types of finite classical polar spaces; Orthogonal polar 
spaces (parabolic quadric $\cQ(2r,q)$,  hyperbolic quadric $\cQ^+(2r-1,q)$, elliptic quadric $\cQ^-(2r-1,q)$); symplectic polar spaces ($\Wa(2r-1,q)$); and Hermitian 
polar spaces ($\Ha(2r,q^2)$, $\Ha(2r-1,q^2)$). See Table~\ref{Tab_Polar} for
polar 
spaces and their ranks and forms $f$. (In  Table~\ref{Tab_Polar}, $f(x_0,x_1)=ax_0^2+bx_0x_1+cx_1^2$ is an irreducible quadratic form in two indeterminates.)
 
\begin{table}[h]
\caption{Classical polar spaces}
\label{Tab_Polar}
\begin{center}
\begin{tabular}{|c|c|c|c|}
\hline 
Polar space& dimension &rank& form\\
\hline \hline 
$\cQ(2r,q)$&$n=2r+1$&$r$&$x_0^2+x_1x_2+\cdots+x_{2r-1}x_{2r}$\\
\hline 
$\cQ^+(2r-1,q)$&$n=2r$&$r$&$x_0x_1+\cdots+x_{2r-2}x_{2r-1}$\\
\hline
$\cQ^-(2r-1,q)$&$n=2r$&$r-1$&$f(x_0,x_1)+x_2x_3+\cdots+x_{2r-2}x_{2r-1}$\\
\hline
$\Wa(2r-1,q)$&$n=2r$&$r$&$x_0y_1+y_0x_1+\cdots+x_{2r-2}y_{2r-1}+x_{2r-1}y_{2r-2}$\\
\hline
$\Ha(2r,q^2)$&$n=2r+1$&$r$&$x_0^{q+1}+\cdots+x_{2r}^{q+1}$\\
\hline
$\Ha(2r-1,q^2)$&$n=2r$&$r$&$x_0^{q+1}+\cdots+x_{2r-1}^{q+1}$\\
\hline
\end{tabular}
\end{center}
\end{table}

Let $S$ be a polar space of rank $r$ over $\F_q$.  
An {\it $m$-ovoid} is a set ${\mathcal M}$ of points of $S$ such that 
every maximal of $S$ meets ${\mathcal M}$ in exactly 
$m$ points. 
For example, the whole point set of $S$ itself is a $\frac{q^r-1}{q-1}$-ovoid. 
For two $m_j$-ovoids ${\mathcal M}_j$, $j=1,2$, if 
${\mathcal M}_2\subseteq {\mathcal M}_1$, then ${\mathcal M}_1\setminus {\mathcal M}_2$ is an $(m_1-m_2)$-ovoid. 
On the other hand, if ${\mathcal M}_1$ and ${\mathcal M}_2$ are disjoint, then  ${\mathcal M}_1\cup {\mathcal M}_2$ is an $(m_1+m_2)$-ovoid. 

For a point $P$ of a polar space $S$, the set $P^\perp$ of points 
of $S$ collinear with $P$ is the intersection of the tangent hyperplane 
at $P$ with $S$. Let ${\mathcal M}$ be an $m$-ovoid of $S$. It is known that  $|P^\perp\cap {\mathcal M}|$ takes  
exactly two values according to $P\in {\mathcal M}$ or not~\cite{BKLP07}. Furthermore, if $S$ is either 
$\Ha(2r,q^2)$, $\cQ^-(2r-1,q)$, or 
$\Wa(2r-1,q)$, the sizes of $H\cap {\mathcal M}$, where $H$ are nontangent hyperplanes, can also be computed exactly. In fact, the following theorem is known.

\begin{theorem}\label{thm:movoid}{\em (\cite[Theorem~11]{BKLP07})}
Let $S$ be one of the polar spaces $\Ha(2r,q^2)$, $\cQ^-(2r-1,q)$, or 
$\Wa(2r-1,q)$ and let ${\mathcal M}$ be an $m$-ovoid in $S$. Then ${\mathcal M}$ is a projective two-intersection set in the ambient projective space of $S$; in other words, letting $
D=\{xy: y \in \F_q^\ast,\langle x\rangle \in {\mathcal M}\}$ and $V$ be the underlying vector space of $S$, the graph 
$\Cay(V,D)$ is an srg with negative Latin square type parameters $(u^2,s(u+1),-u+s^2+3s,s^2+s)$, where 
$(u,s)=(q^{2r+1},m(q^2-1))$, $(q^r,m(q-1))$, or $(q^r,m(q-1))$ according as  $S=\Ha(2r,q^2)$, $\cQ^-(2r-1,q)$, or 
$\Wa(2r-1,q)$, respectively. 
\end{theorem}

\begin{remark}
{\em 
\item[(1)] A partial converse to the above theorem holds. That is, if ${\mathcal M}$ is a projective two-intersection set in the ambient projective space of $S$, and ${\mathcal M}$ satisfies certain conditions, then ${\mathcal M}$ is an $m$-ovoid in $S$. We refer the reader to \cite{BLMX18} for the precise statement of the partial converse. This partial converse provides an approach to constructing $m$-ovoids in the polar spaces mentioned in Theorem~\ref{thm:movoid}.
\item[(2)] A $(q+1)/2$-ovoid in 
$\cQ^-(5,q)$ can be interpreted as a set of lines in $\Ha(3,q^2)$ containing 
exactly half of the lines on every point via the duality of generalized 
quadrangles. Such a set of lines in $\Ha(3,q^2)$ is called a {\it hemisystem}, which was first studied by Segre~\cite{Se65}. Constructions of hemisystems can be found in \cite{BGR10,BLMX18,CP05,KNS17}. }
 %
\end{remark}

To obtain a similar theorem for srgs of Latin square type, we need to introduce the concept of $i$-tight sets.  Let $S$ be a polar space of 
rank $r \ge 2$ over $\F_q$. An {\it $i$-tight set} is a set ${\mathcal M}$ of points of $S$ 
such that 
\[
|P^\perp\cap {\mathcal M}|=
\begin{cases}
i\frac{q^{r-1}-1}{q-1}+q^{r-1},& \text{ if $P\in {\mathcal M}$}, \\
i\frac{q^{r-1}-1}{q-1}, & \text{ otherwise.}
\end{cases}
\]
For example, 
each maximal is a $1$-tight set. 
In \cite{BKLP07}, it was shown that if a set ${\mathcal M}$ of points in a polar space $S$ meets $P^\perp$ in exactly two different sizes according to $P\in {\mathcal M}$ or not, then ${\mathcal M}$ is either an $m$-ovoid or an $i$-tight set for some $m$ or $i$. Similarly to the situation with $m$-ovoids, the following basic properties hold. 
For two $i_j$-tight sets ${\mathcal M}_j$ in $S$, $j=1,2$, if 
${\mathcal M}_2\subseteq {\mathcal M}_1$, then ${\mathcal M}_1\setminus {\mathcal M}_2$ is an $(i_1-i_2)$-tight set. 
On the other hand, if ${\mathcal M}_1$ and ${\mathcal M}_2$ are disjoint, then  ${\mathcal M}_1\cup {\mathcal M}_2$ is an $(i_1+i_2)$-tight set. 
Furthermore, if $S$ is either 
$\Ha(2r-1,q^2)$, $\cQ^+(2r-1,q)$, or 
$\Wa(2r-1,q)$, the size of $H\cap {\mathcal M}$ for nontangent hyperplanes can be also computed exactly. In fact,  the following theorem is known. 
\begin{theorem}\label{thm:itight}{\em (\cite[Theorem~12]{BKLP07})}
Let $S$ be one of the polar spaces $\Ha(2r-1,q^2)$, $\cQ^+(2r-1,q)$, or 
$\Wa(2r-1,q)$ and let ${\mathcal M}$ be an $i$-tight set in $S$. 
Then ${\mathcal M}$ is a projective two-intersection set in the ambient projective space of $S$; In other words, 
letting $D:=\{xy: y \in \F_q^\ast,\langle x\rangle \in {\mathcal M}\}$ and $V$ be the underlying vector space of $S$, 
the graph $\Cay(V,D)$ is an srg  with Latin square type parameters $(u^2,s(u-1),u+s^2-3s,s^2-s)$, where 
$(u,s)=(q^{2r},i)$, $(q^r,i)$, or $(q^r,i)$ according as  $S=\Ha(2r-1,q^2)$, $\cQ^+(2r-1,q)$, or 
$\Wa(2r-1,q)$, respectively. 
\end{theorem}

\begin{remark}
{\em
\item[(1)] Again, a partial converse to the above theorem holds. For the detailed statement, see \cite{FMX152}. This partial converse provides an approach to constructing Cameron-Liebler lines classes in $\PG(3,q)$ (see definition below).
\item[(2)] A tight set in $\cQ^+(5,q)$
can be interpreted as a set ${\mathcal L}$ of lines in $\PG(3,q)$ such that the size of ${\mathcal L}\cap S$ is constant for all spread $S$ via the Klein 
correspondence.  Such a set of lines in $\PG(3,q)$ is called a {\it Cameron-Liebler line class}, which was first studied by Cameron and Liebler~\cite{CL82}. Constructions of Cameron-Liebler line classes can be found in  \cite{BD99,DDMR16,FMX152,GMP18,Rod12}. }
\end{remark}

Known results on $m$-ovoids and $i$-tight sets are surveyed in \cite{BKLP07,BLP09}. 
See \cite{CCEM08,CP14,CP17,De17,NS17,Me16} for recent constructions of $m$-ovoids and $i$-tight sets.

\providecommand{\bysame}{\leavevmode\hbox to3em{\hrulefill}\thinspace}
\providecommand{\MR}{\relax\ifhmode\unskip\space\fi MR }
\providecommand{\MRhref}[2]{%
  \href{http://www.ams.org/mathscinet-getitem?mr=#1}{#2}
}
\providecommand{\href}[2]{#2}


\end{document}